\documentclass{amsart}
\usepackage{amsmath,amsthm}
\usepackage{amsfonts,amssymb}

\usepackage{enumerate}
\usepackage{multirow}

\newtheorem{thm}{Theorem}[section]

\newtheorem{prop}[thm]{Proposition}

\theoremstyle{remark}

 \def\a{{\alpha}} 
 \def\b{{\beta}}
 \def\g{{\gamma}}

 \def\l{{\lambda}}

 \def\CL{{\mathcal L}}

 \def\NN{{\mathbb N}}

\newif\ifpdf
\ifx\pdfoutput\undefined
  \pdffalse
\else
  \pdftrue
\fi

\ifpdf
  \usepackage[pdftex]{graphicx}
  \DeclareGraphicsExtensions{.pdf,.jpg,.png}
\else
  \usepackage{graphicx}
\fi

\begin{document}
 
\title[Gauss-Lobatto integration on the triangle]
{On Gauss-Lobatto integration on the triangle}
\author{Yuan Xu}
\address{Department of Mathematics\\ University of Oregon\\
    Eugene, Oregon 97403-1222.}\email{yuan@math.uoregon.edu}

\date{\today}

\begin{abstract}
A recent result  in [2] on the non-existence of Gauss-Lobatto cubature rules on the triangle 
is strengthened by establishing a lower bound for the number of nodes of such rules. 
A method of constructing Lobatto type cubature rules on the triangle is given and used to 
construct several examples.   
\end{abstract}

\maketitle

\section{Introduction}
\setcounter{equation}{0}

Recently in \cite{H}, motivated by $hp$-finite element method, Gauss-Lobatto 
cubature rule on the triangle
$$
  \triangle := \{(x,y): 0 \le x, y, x+y \le 1\}
$$
is studied, which requires $(n+1)(n+2)/2$ nodes, with $(n-1)(n-2)/2$ nodes in 
the interior, $n-1$ nodes in each side and $1$ point at each vertex, and is capable 
of exactly integrating  polynomials of degree $2n-1$. We call such a rule strict
Gauss-Lobatto. The main result of \cite{H} shows that such rules do not exist. 
In this note we consider cubature rules that have nodes in both interior and 
boundary of the triangle, and establish a lower bound for the number of nodes, from 
which the non-existence of the strict Gauss-Lobatto rule follows immediately. We also 
study the structure of rules that attain our lower bound and give a method for 
constructing cubature rules of degree $2n-1$ with $n-1$ nodes on each side and 
1 node at each vertex. The development is based on the observation that cubature
rules with nodes on the boundary can be constructed, by restricting to the class
of bubble functions (functions that vanish on the boundary), from 
cubature rules with nodes on the interior. This leads to a bootstrapping scheme for 
transforming some cubature rules with interior points into higher order rules with a 
specific number of nodes on the boundary. Several examples are constructed to 
illustrate the algorithm.

\section{Results}
\setcounter{equation}{0}

For $n \in \NN_0$, let $\Pi_n^2$ denote the space of 
polynomials of (total) degree $n$ in two variables. It is known that 
$$
         \dim \Pi_n^2 = \binom{n+2}{2} = \frac{(n+1)(n+2)}{2}. 
$$
Let $W(x,y)$ be a non-negative weight function on the triangle $\triangle$ 
with finite moments. A cubature rule of precision $s$ with respect to $W$ is a finite
sum that satisfies
\begin{equation} \label{cubaa}
   \int_\triangle f(x,y) W(x,y) dxdy = \sum_{k=1}^N \l_k f(x_k,y_k), \qquad \forall f\in \Pi_s^2. 
\end{equation}
We choose $W$ to be the Jacobi weight $W_{\a,\b,\g} (x,y) = x^\a y ^\b (1-x-y)^\g$ for 
$\a,\b,\g > -1$. The case $\a = \b =\g =0$ corresponds to the constant weight. These 
weight functions are often considered together with the orthogonal polynomials, called 
Jacobi polynomials on the triangle, that are orthogonal with respect to them;
see, for example, \cite[p. 86]{DX}, and \cite{BS, BSX} in connection with cubature rules. 

For the Jacobi weight, it is known that the number of nodes $N$ for \eqref{cubaa} satisfies
\begin{equation}\label{lwbd}
   N \ge \begin{cases}    \frac{n(n+1)}{2} +  \lfloor \frac{n}{2}\rfloor & \hbox{if $s = 2n-1$},\\
        \frac{n(n+1)}{2} & \hbox{if $s = 2n-2$}. \end{cases}  
\end{equation}
This lower bound is classical for $s = 2n-2$ (\cite{St}) and given in \cite{BS} for 
$s = 2n-1$, which agrees with M\"oller's lower bound for centrally symmetric weight
functions \cite{M} as well as \cite{DX}. A cubature rule that attains the lower bound is naturally 
minimal, meaning that it has the smallest number of nodes among all cubature rules of the 
same degree. The lower bound, however, is most likely not sharp; that is, a minimal cubature 
could require more points than what the lower bound indicates. Minimal cubature rules are
sometimes called Gaussian cubature rules. Their construction is closely related to
orthogonal polynomials of several variables. For discussion along this line, see \cite{CMS, X}
and references therein. 

The cubature rules that we consider are of precision $s$ and are of the form
\begin{align} \label{cuba}
   \int_\triangle f(x,y) W_{\a,\b,\g}(x,y) dxdy  & = \sum_{k=1}^{N_0} \l_{k,0} f(x_{k,0},y_{k,0})   
    + \sum_{k=1}^{N_1} \l_{k,1} f(x_{k,1},0) \\ 
    & + \sum_{k=1}^{N_2} \l_{k,2} f(0,y_{k,2}) +    \sum_{k=1}^{N_3} \l_{k,3} f(x_{k,3}, 1-x_{k,3}) 
      \notag \\ 
    & + \mu_0 f(0,0) + \mu_1 f(1,0) + \mu_2 f(0,1),   \notag   
\end{align}
where $(x_{k,0},y_{k,0})$ are distinct points in the interior of $\triangle$, $(x_{k,1}, 0)$, $(0,y_{k,2})$, 
and $(x_{k,3},1-x_{k,3})$ are distinct points on the side $y=0$, $x=0$, and $x+y=1$ 
(but not on the corners) of $\triangle$, respectively, $\lambda_{k,j} > 0$ and $\mu_i > 0$. 
Such a cubature has 
$$
       N : = N_0 + N_1 +N_2+N_3 + 3
$$
nodes. The main result in \cite{H} states that such a cubature rule does not exist if 
$s = 2n-1$ and 
$$
 N_0 = \frac{(n-2)(n-1)}{2}, \quad N_1 = N_2 =N_3  = n-1, 
$$
which has a total number of nodes $N = (n+2)(n+1)/2$. This follows as an immediate corollary of the following theorem.

\begin{thm}
If a cubature rule of the form \eqref{cuba} exists with precision $s = 2n-1$ or $s =2n$, then 
\begin{align}
   N_0 & \ge \begin{cases} \frac{n(n-1)}{2} & \hbox{if $s = 2n-1$}, \\
          \frac{n (n-1)}{2} +  \lfloor \frac{n-1}{2}\rfloor & \hbox{if $s = 2n$},\end{cases} \label{lwbd2} \\
   N_0 + N_i & \ge \begin{cases}  \frac{n(n-1)}{2}+  \lfloor \frac{n-1}{2}\rfloor & \hbox{if $s = 2n-1$}, \\
        \frac{n (n+1)}{2}  & \hbox{if $s = 2n$},\end{cases}  \qquad i =1,2,3. \label{lwbd2b}
\end{align}
\end{thm}
  
\begin{proof}
The cubature \eqref{cuba} exactly integrates degree $s$ polynomials of the form 
$f(x,y) = xy(1-x-y) g(x,y)$ if 
\begin{equation} \label{cuba2}
  \int_\triangle g(x,y)  W_{\a+1,\b+1,\g+1} (x,y)dxdy = \sum_{k=1}^{N_0}
      \l_{k,0}^* g(x_{k,0}, y_{k,0}), \quad  \forall g \in \Pi_{s-3}^2, 
\end{equation}
where $\l_{k,0}^* = \l_{k,0} x_{k,0}y_{k,0} (1-x_{k,0} - y_{k,0})$, which is a cubature rule
of precision $s-3$ for the weight function $W_{\a+1,\b+1,\g+1}$ so that, by \eqref{lwbd}, 
$N_0$ has to satisfy the lower bound in the inequality of \eqref{lwbd2}. 
On the other hand, 
the cubature (2.3) exactly integrates degree $s$ polynomials of the form
$f(x,y)=x(1-x-y)g(x,y)$ if $\forall g \in \Pi_{s-3}^2$,
\[
\int_{\Delta} g(x,y) W_{\alpha+1,\beta,\gamma+1}(x,y)dxdy=
\sum_{k=1}^{N_0} \tilde{\lambda}_{k,0} g(x_{k,0}, y_{k,0})
+
\sum_{k=1}^{N_1} \lambda_{k,1}^* g(x_{k,1},0)
\]
where $\tilde{\lambda}_{k,0} = \lambda_{k,0} x_{k,0} (1-x_{k,0}-y_{k,0})$
and $\lambda_{k,1}^* = \lambda_{k,1} x_{k,1} (1-x_{k,1})$, which is a cubature rule of precision $s-2$
for the weight function $W_{\a+1,\b,\g+1}(x,y)$, so that $N_0+N_1$ satisfies the lower bound
in \eqref{lwbd},  which gives the inequality of \eqref{lwbd2b} for $i=1$. Similarly, we 
can derive lower bound for $N_0+N_2$ and $N_0+N_3$. 
\end{proof}

One naturally asks if there is any cubature rule that attains the lower bound in the theorem. 
For $s = 2n-1$, this asks if there is a cubature rule of precision $2n-1$ with 
\begin{align}  \label{min} 
 N_0  =  \frac{n(n-1)}{2} \quad \hbox{and} \quad  N_i = \left \lfloor \frac{n-1}{2} \right \rfloor,  \quad i =1,2,3. 
\end{align}
We expect that the answer is negative. A heuristic argument can be given as follows: Assume that
$N_0 =\frac{n(n-1)}{2}$. Then the proof of the theorem shows that \eqref{cuba2} is a cubature of 
degree $2n-4$ with $n(n-1)/2$ nodes, which is known to exist only for small $n$. Assume that 
it does exist. We define a linear functional $\CL_1$, acting on polynomials of one variable, by 
\begin{equation}\label{L1}
 \CL_1 g : = \int_\triangle g(x)  W_{\a+1,\b,\g+1}(x,y)dxdy - 
    \sum_{k=1}^{N_0} \l_{k,0} x_{k,0} (1-x_{k,0} - y_{k,0}) g(x_{k,0}),
\end{equation}
where $x_{k,0}$, $y_{k,0}$ and $\l_{k,0}$ are as in \eqref{cuba}. Applying \eqref{cuba} on 
polynomials of the form $f(x,y) = g(x) x (1-x-y)$ shows that 
\begin{equation}\label{gausian}
  \CL_1 g = \sum_{k=1}^{N_1} \l_{k,1}^* g(x_{k,1}), \qquad  \forall g \in \Pi_{2n-3}, 
\end{equation}
where $\l_{k,1}^* = \l_{k,1} x_{k,1}(1-x_{k,1})$. 
The functional $\CL_1$ defined in equation \eqref{L1} defines a bilinear form $[p,q]_1= \CL_1(pq)$ 
that could be indefinite ($[q,q]_1$ is not necessarily positive). If the bilinear form were positive definite 
on $\Pi_{2n-3}$, then \eqref{gausian} could be regarded as a quadrature rule of $N_1$ nodes 
and of degree $2n-3$ for $\CL_1$ and, consequently, $N_1 \ge n-1$ by the standard result in Gaussian 
quadrature rule, which is stronger than the second equation of \eqref{min}. Thus, in order for 
\eqref{min} to hold,  we would need $\CL_1$ to be indefinite on $\Pi_{2n-3}$, that is, $\CL (q^2) = 0$
for some nonzero $q \in \Pi_{n-1}$, and we would need to require that $\CL_1$ has a quadrature rule of 
degree $2n-3$ with $N_1 = \lfloor \frac{n-1}{2} \rfloor$ nodes. A simple count of variables (the 
nodes and weights of the quadrature rule) and restraints (the polynomials that need to be exactly
integrated) shows that this is unlikely to happen, although it still might as the equations are nonlinear. 
For a linear functional that defines an indefinite 
bilinear form, the theory of Gaussian quadrature rule breaks down since orthogonal polynomials 
may not exist and, even they do, they may not have real or simple zeros. In particular, we do not 
have a lower bound for the number of nodes of a quadrature rule for such a linear functional.  

The above argument indicates that if we want the cubature rule \eqref{cuba} to have the smallest
number of interior points, then the best that we can hope for will be a cubature rule of degree $2n-1$
that satisfies 
\begin{equation}\label{min2} 
      N_0  =  \frac{n(n-1)}{2} , \quad \hbox{and} \quad N_i \ge n -1, \qquad 1 \le i \le 3. 
\end{equation}

We shall call a cubature rule that attains the lower bound in \eqref{min2} {\it Gauss-Lobatto 
cubature rule}. The proof of the lower bound indicates how such a cubature rule can be constructed.  
Let us also define linear functionals 
\begin{align*}
 \CL_2 g &: = \int_\triangle g(y)  W_{\a,\b+1,\g+1}(x,y)dxdy - 
    \sum_{k=1}^{N_0} \l_{k,0} y_{k,0} (1-x_{k,0} - y_{k,0}) g(x_{k,0}), \\
 \CL_3 g & : = \int_\triangle g(x)  W_{\a+1,\b+1,\g}(x,y)dxdy - 
    \sum_{k=1}^{N_0} \l_{k,0} x_{k,0}  y_{k,0} g(x_{k,0}).
\end{align*}
Then we can summarize the method of construction as follows. 

\medskip\noindent
{\bf Algorithm.} We can follow the following procedure to construct a Gauss-Lobatto cubature rule
of form \eqref{cuba}: 

{\it Step 1.} Construct a cubature rule of degree $s-3$ for $W_{\a+1,\b+1,\g+1}$ in the form of
\eqref{cuba2} with all nodes in the interior of $\triangle$, and define 
\begin{equation} \label{lambda0}
  \l_{k,0} = \frac{ \l_{k,0}^*}{(x_{k,0}y_{k,0}(1-x_{k,0}-y_{k,0}))}.
\end{equation}
{\it Step 2.} Construct Gaussian quadrature rules \eqref{gausian} and 
\begin{equation*}
  \CL_2 g = \sum_{k=1}^{N_2} \l_{k,2}^* g(y_{k,2}), \quad \CL_3 g = \sum_{k=1}^{N_3} \l_{k,3}^* g(x_{k,3}), 
  \qquad \forall g \in \Pi_{2n-3}, 
\end{equation*}
with respect to the linear functional $\CL_1,\CL_2,\CL_3$, which gives the nodes of the cubature rule
\eqref{cuba} on the boundary of the triangle, and define 
\begin{equation} \label{lambda*}
\l_{k,1} =  \frac{ \l_{k,1}^*}{x_{k,1}(1-x_{k,1})}, \quad \l_{k,2} =  \frac{\l_{k,2}^*}{ y_{k,2}(1-y_{k,2})}, \quad
 \l_{k,3} = \frac{\l_{k,3}^*}{x_{k,3}(1-x_{k,3})}. 
\end{equation}
 
{\it Step 3.} Finally, the weight $\mu_0$, $\mu_1$, and $\mu_2$ are determined by setting
$f(x,y) = 1, x, y$ in \eqref{cuba} and solve the resulted linear system of equations. 
 \qed

\begin{prop} \label{prop1}
Assume, in the above algorithm, that the linear functional $\CL_1,\CL_2$ and $\CL_3$ are positive definite 
on $\Pi_{2n-3}$ and all $x_{k,1}, y_{k,2}, x_{k,3}$ are inside $(0,1)$. Then the algorithm produces a cubature 
rule of degree $2n-1$ in the form \eqref{cuba}. In particular, if $N_0 = \frac{n(n-1)}{2}$, then the cubature rule 
is Gauss-Lobatto.
\end{prop}

\begin{proof}
We need to show that the cubature rule constructed by the algorithm holds for all $f \in \Pi_{2n-1}^2$, 
that is, \eqref{cuba} holds for all $f \in \Pi_s^2$ with $s = 2n-1$. A moment of reflection shows that, with 
$z = 1-x-y$, $\Pi_s^2$ can be decomposed into a direct sum.
\begin{equation}\label{decom}
 \Pi_s^2 =x y z \Pi_{s-3}^2 \oplus xz\Pi_{s-2}[x]\oplus yz\Pi_{s-2}[y]\oplus x y\Pi_{s-2}[x]\oplus \Pi_1^2, 
\end{equation}
where $\Pi_{s-2}[x]$ and $\Pi_{s-2}[y]$ denote the space of polynomials of one variable in $x$-variable
and $y$-variable, respectively. If $f \in x y z \Pi_{s-3}^2$, then \eqref{cuba} reduces to \eqref{cuba2},
which holds by our construction. If $f \in xz\Pi_{s-2}[x]$, then \eqref{cuba} reduces to \eqref{gausian} for
$\CL_1$, which holds by our construction. The same holds for $f \in yz\Pi_{s-2}[y]$ and $f \in  x y\Pi_{s-2}[x]$, 
whereas for $f \in \Pi_1^2$, the cubature is verified by Step 3. Thus, by \eqref{decom}, the cubature
rule holds for all $f\in \Pi_{2n-1}$. 
\end{proof}

It should be pointed out that, as long as we have a cubature rule of degree $2n-4$ with all nodes in the 
interior of $\triangle$ in the Step 1, regardless if it is a minimal one, the Step 2 and Step 3 could be 
carried out and Proposition \ref{prop1} applies. Thus, the algorithm can be used to construct cubature 
rules of degree $2n-1$ in the form \eqref{cuba} with $N_1 = N_2 = N_3 = n-1$. 

Let us comment on how the steps in the algorithm can be realized. 
 
For Step 1, a cubature with the specification can be constructed by solving moment equations, 
that is, solving the system of equations formed by setting $g(x,y) = x^iy^j$ for $i+j \le s-3$ in 
\eqref{cuba2} for $x_{k,0}$, $y_{k,0}$ and $\l_{k,0}^*$. There have been a number of papers based
on this method, see e.g. \cite{ZCL} for the latest result and further references. Another approach is to use
a characterization of the cubature rules that attain lower bound in \eqref{lwbd2}, which is given in
terms of common zeros of certain orthogonal polynomials and can be used to find cubature rules of 
lower order, see \cite{M, X}. Not all cubature rules obtained via either methods work for our purpose,
since we require that all nodes are inside the domain. 

For Step 2, one needs to check that $\CL_1, \CL_2,\CL_3$ are positive linear functionals on $\Pi_{2n-3}$.
Once they are, the standard procedure of constructing Gaussian quadrature rules applies. In 
particular, we can apply the standard algorithm to generate a sequence of orthogonal polynomials 
up to degree $n-1$ with respect to $\CL_i$ inductively; the nodes of the Guassian quadrature rule 
of degree $2n-3$ for $\CL_i$ are the zeros of the orthogonal polynomial of degree $n-1$ with 
respect to $\CL_i$.

\medskip

In the following we give several examples of Gauss-Lobatto cubature rules of degree $5$ and $7$ for 
the unit weight function. The minimal cubature rules for the degree 5 and 7 have nodes 7 and 12, 
respectively (\cite{CR}). For Gauss-Lobatto rules, the number of nodes are necessarily larger, as 
seen in \eqref{min}. 

\medskip\noindent
{\bf Example 1.} There is a cubature formula of degree $5$ with 12 nodes, 3 interior, 2 on each 
side and 1 at each corner of the triangle. To illustrate our procedure, we shall present the nodes 
and weights in steps. 

The three interior points and weights are given in the first table. While the nodes are those of a 
cubature rule of degree 2 for the weight function $W_{1,1,1}(x,y) = xy(1-x-y)$, the weights are
relates to those of the latter cubature rule by \eqref{lambda0}. These nodes are common zeros 
of quasi-orthogonal polynomials (see \cite{X2} for definition) of degree $2$ which were found by
solving the nonlinear system of equations in Theorem 4.1 of \cite{X2}. 

\medskip

\hskip .5cm \vbox{\tabskip=0pt \offinterlineskip
\def\tablerule{\noalign{\hrule}}
\halign to 310pt{\strut#& \vrule#\tabskip=1em 
   \hfil& \hfil#\hfil& \vrule#& \hfil #\hfil& \vrule#& \hfil#\hfil& \vrule#&
   \hfil#& \vrule#\tabskip=0pt \cr\tablerule
&& \omit\hidewidth $x_{k,0}$ \hidewidth&&
  \omit\hidewidth  $y_{k,0}$ \hidewidth&&
  \omit\hidewidth  $\lambda_{k,0}$ \hidewidth&\cr\tablerule
 &&0.15881702219143 && 0.19201873632215 && 0.101342396527698 &\cr\tablerule
 &&0.56219234596964 && 0.19201873632215 && 0.117181247909596  &\cr\tablerule
 &&0.22100936816107 && 0.55798126367785 && 0.118066904793533 &\cr\tablerule
}}  

\noindent where the nodes are the common zeros of following three polynomials of degree $2$:
\begin{align*}
 & 19 - \sqrt{105}+ x (-91 + \sqrt{105} + 112 x) + 2 (-7 + \sqrt{105}) y, \\ 
 & 49 - 3 \sqrt{105} - 175 x + 5 \sqrt{105} x + 112 x^2 - 84 y +  4 \sqrt{105} y + 224 x y, \\ 
 &154 - 6 \sqrt{105} - 301 x + 11 \sqrt{105} x + 112 x^2 - 609 y +  7 \sqrt{105} y + 560 x y + 560 y^2.
\end{align*}

\medskip

\noindent
The nodes on the edges of the triangle, but not on the corners, are the nodes of Gaussian 
quadrature rules for $\CL_1, \CL_2, \CL_3$, respectively, where the weights are related to
those of Gaussian quadrature rules by \eqref{lambda*}. They are given in the following three
tables:

\medskip
 
\hskip .5cm \vbox{\tabskip=0pt \offinterlineskip
\def\tablerule{\noalign{\hrule}}
\halign to 199pt{\strut#& \vrule#\tabskip=1em 
   \hfil& \hfil#\hfil& \vrule#& \hfil #\hfil& \vrule#& \hfil#\hfil& \vrule#&
   \hfil#& \vrule#\tabskip=0pt \cr\tablerule
&& \omit\hidewidth $x_{k,1}$ \hidewidth&&  \omit\hidewidth  $\lambda_{k,1}$ \hidewidth&\cr\tablerule
 && 0.3931870086016  && 0.02991955921794   &\cr\tablerule
 && 0.8595419130359 & & 0.01756588222187  &\cr\tablerule
}}

\noindent 
where the nodes are zeros of the orthogonal polynomials  $p_1$ of degree 2 for $\CL_1$,
$$
 p_1(x) = x^2 + \frac{1}{448} \left(-469 - 9 \sqrt{105}\right) x  + \frac{1}{4480}\left(889 + 61 \sqrt{105}\right);
$$

\medskip

\hskip .5cm \vbox{\tabskip=0pt \offinterlineskip
\def\tablerule{\noalign{\hrule}}
\halign to 199pt{\strut#& \vrule#\tabskip=1em 
   \hfil& \hfil#\hfil& \vrule#& \hfil #\hfil& \vrule#& \hfil#\hfil& \vrule#&
   \hfil#& \vrule#\tabskip=0pt \cr\tablerule
&& \omit\hidewidth $y_{k,2}$ \hidewidth&&  \omit\hidewidth  $\lambda_{k,2}$ \hidewidth&\cr\tablerule
 && 0.4305843026985 && 0.02290932968619   &\cr\tablerule
 && 0.7924406473476 && 0.02022650113138   &\cr\tablerule
}}

\noindent 
where the nodes are zeros of the orthogonal polynomials $p_2$ of degree 2 for $\CL_2$,
$$
 p_2(x) =  x^2 + \frac{3}{46} \left(-29 + \sqrt{105}\right) x  + \frac{3}{644} \left (63 + \sqrt{105}\right);
$$

\medskip

\hskip .5cm \vbox{\tabskip=0pt \offinterlineskip
\def\tablerule{\noalign{\hrule}}
\halign to 199pt{\strut#& \vrule#\tabskip=1em 
   \hfil& \hfil#\hfil& \vrule#& \hfil #\hfil& \vrule#& \hfil#\hfil& \vrule#&
   \hfil#& \vrule#\tabskip=0pt \cr\tablerule
&& \omit\hidewidth $x_{k,3}$ \hidewidth&&  \omit\hidewidth  $\lambda_{k,3}$ \hidewidth&\cr\tablerule
 && 0.2629899118578 &&   0.02514330117112  &\cr\tablerule
 && 0.7030163143652  &&  0.03109870484395  &\cr\tablerule
}}

\noindent 
where the nodes are zeros of the orthogonal polynomials $p_3$ of degree 2 for $\CL_3$,
$$
 p_3(x) =  x^2 + \frac{1}{10843} \left(-10997 + 51 \sqrt{105} \right) x + \frac{1}{3098}\left(665 - 9 \sqrt{105}\right).
$$
Finally, the weights $\mu_1,\mu_2,\mu_3$ are given by 

\medskip

\hskip .5cm \vbox{\tabskip=0pt \offinterlineskip
\def\tablerule{\noalign{\hrule}}
\halign to 300pt{\strut#& \vrule#\tabskip=1em 
   \hfil& \hfil#\hfil& \vrule#& \hfil #\hfil& \vrule#& \hfil#\hfil& \vrule#&
   \hfil#& \vrule#\tabskip=0pt \cr\tablerule && \omit\hidewidth $\mu_1$ \hidewidth&&
  \omit\hidewidth  $\mu_2$ \hidewidth&&  \omit\hidewidth  $\mu_3$ \hidewidth&\cr\tablerule
 && 0.0081170837035 &&0.00326155091683 && 0.00516753787639 &\cr\tablerule
 }}  \qed
 
The nodes of this Lobatto cubature rule are depicted in Figure 1.

\begin{figure}[ht] 
\centerline{\includegraphics[width = 5cm]{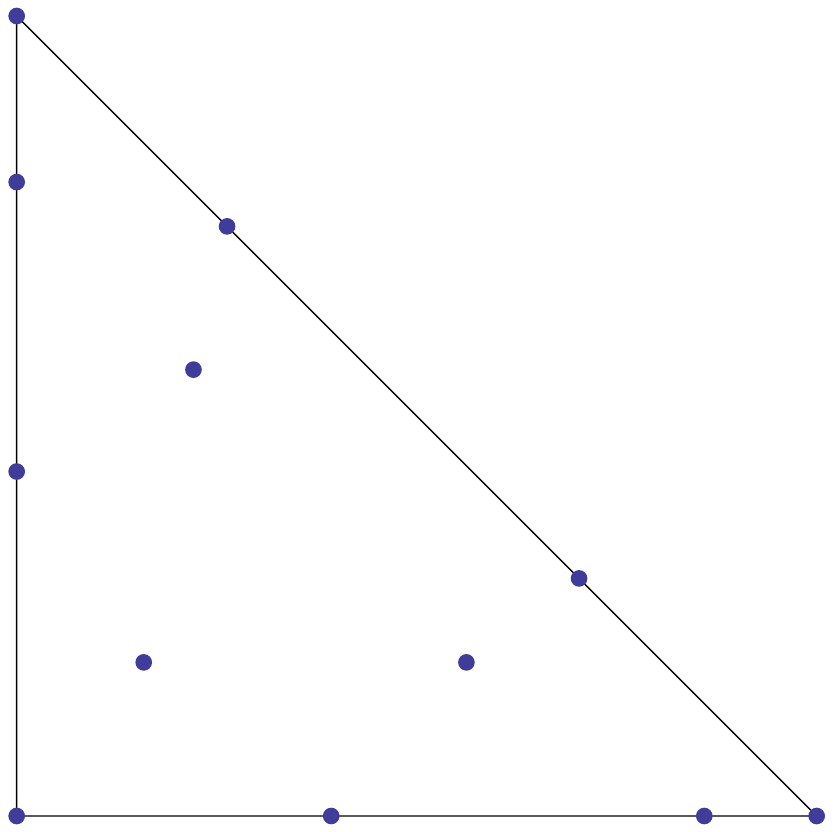} }
\caption{Nodes of Gauss-Lobatto cubature rules of degree 5} \label{Fig1}
\end{figure} 
 
Our next two examples are Gauss-Lobatto cubature rules that are symmetric in the sense
that the nodes are invariant under the symmetric group of $\triangle$, or the permutation of 
$(x,y,1-x-y)$. A symmetric Gauss-Lobatto rule for the constant weight takes the form
\begin{align} \label{symm}
  \int_\triangle f(x,y) dxdy  & =  \sum_{k=1}^{M_0} A_{k} \left[ f(u_{k,0},v_{k,0}) + f(v_{k,0},w_{k,0}) 
         + f(w_{k,0},u_{k,0})\right] \\  
         & + \sum_{k=1}^{M_1} B_{k} \left[ f(u_{k,1},0) + f(0,1-u_{k,1}) + f(1-u_{k,1}, u_k) \right]  \notag \\  
         & + C  \left[ f(0,0) + f(1,0) + f(0,1) \right], \notag
\end{align}
where $w_{k,0} = 1-u_{k,0} - v_{k,0}$. 

Both examples are constructed by following the steps in the algorithm. The symmetry makes the construction 
much easier, since we only need to consider polynomials that are symmetric under the symmetric group of 
$\triangle$. In particular, the cubature rule for $W_{1,1,1}(x,y)$ in Step 1 can be found by solving the reduced 
moment equations of symmetric polynomials. We shall skip details and only list the nodes and weights of 
these two cubature rules as formulated in \eqref{symm}.  Their nodes are depicted in Figure 2. 

\medskip\noindent
{\bf Example 2.} Symmetric Gauss-Lobatto cubature rules of degree $5$ with 12 nodes. This formula 
is in the form of \eqref{symm} with $M_0=1$ and $M_1 =2$. The nodes and weights are given below: 
\begin{align*}
  & u_{1,0}  =  v_{1,0} =  \tfrac{1}{21} (7 - \sqrt{7}), \quad  A_1=  \tfrac{7}{720} \big(14 - \sqrt{7}\big),  \\
  & u_{1,1} =  \tfrac{1}{42} \Big(21 - \sqrt{21 \big( 4 \sqrt{7}-7\big)}\Big),  \quad 
     u_{1,2} = \tfrac{1}{42} \Big(21 + \sqrt{21 \big( 4 \sqrt{7}-7\big)}\Big),  \\
  & B_1=   B_2 = \tfrac{1}{720} \big(7 + 4 \sqrt{7}\big), 
    \quad C = \tfrac{1}{720} \big(8 - \sqrt{7}\big).  \end{align*}  \qed

\medskip\noindent
{\bf Example 3.} Symmetric Gauss-Lobatto cubature rules of degree $7$ with 18 nodes. This formula 
is in the form of \eqref{symm} with $M_0=2$ and $M_1 =3$. The nodes and weights are given below: 
\begin{align*}
 & u_{1,0}  =  v_{1,0} = \tfrac{1}{18} (5 - \sqrt{7}),   \qquad  u_{2,0}  =  v_{2,0} = \tfrac{1}{18} (5 + \sqrt{7}), \\
 & A_1=  \tfrac{1}{17640} \big(1141 -94 \sqrt{7}\big), \quad A_2 =  \tfrac{1}{17640} \big(1141 + 94 \sqrt{7}\big), \\
 & u_{1,1} = \tfrac{1}{6} \big(3 - \sqrt{3}\big), \quad u_{2,1} = \tfrac{1}{2}, \quad u_{3,1} = \tfrac{1}{6} \big(3 + \sqrt{3}\big),\\
 & B_1 = \tfrac{3}{280}, \quad B_2 = \tfrac{4}{315}, \quad B_3= \tfrac{3}{280},   \quad C = \tfrac{1}{315}.
\end{align*} \qed

\begin{figure}[ht] 
\centerline{\includegraphics[width = 5cm]{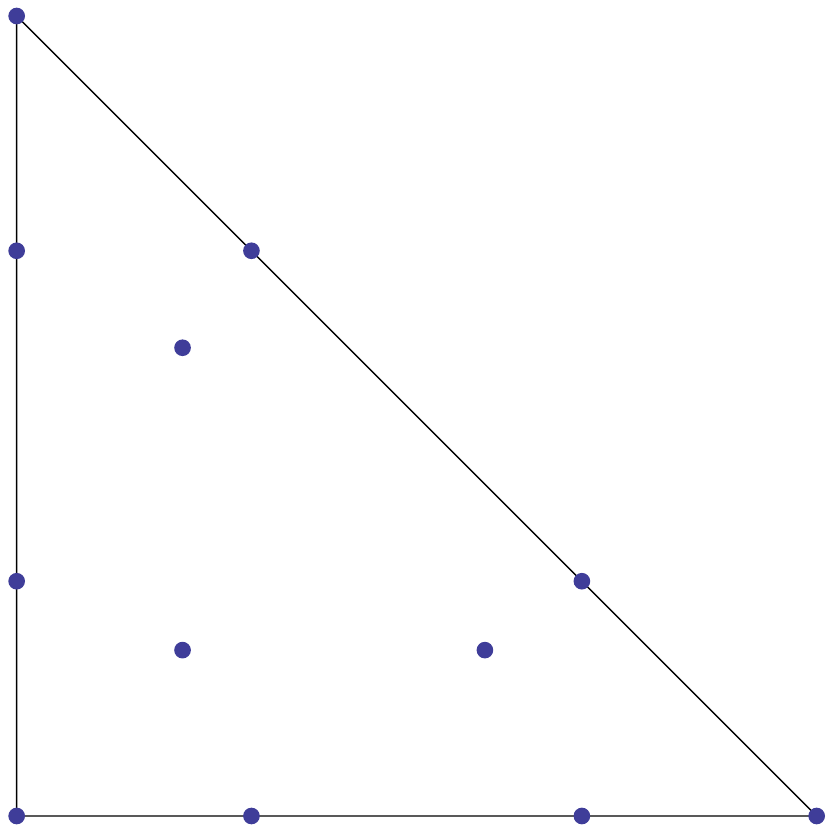} \,\, \,\,
   \includegraphics[width = 5cm]{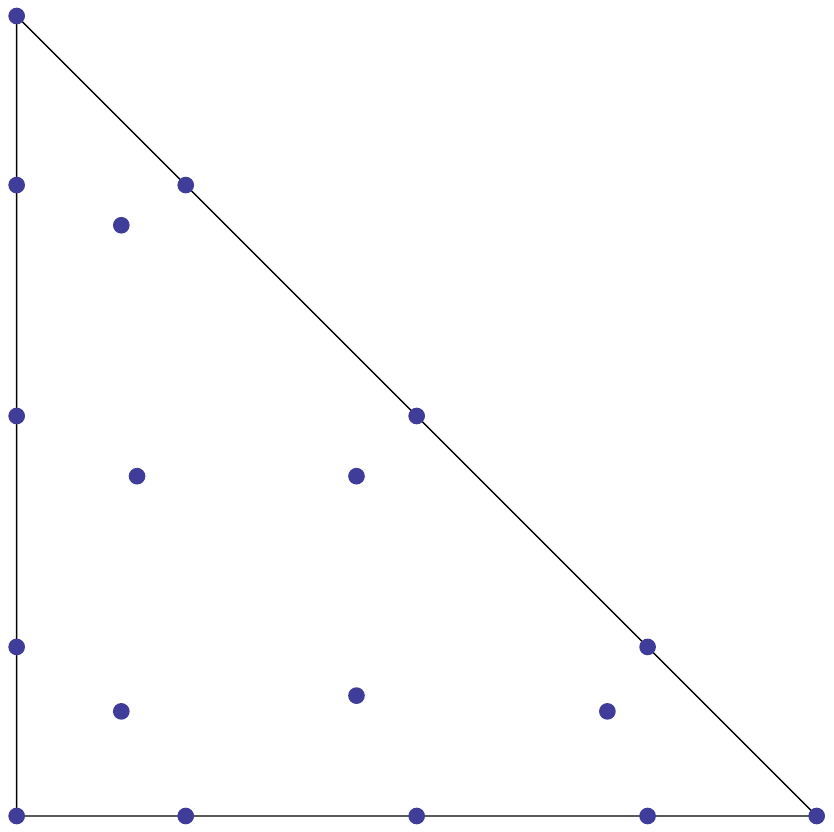} }
\caption{Nodes of symmetric Gauss-Lobatto cubature rules of degree 5 and 7} \label{Fig2}
\end{figure} 

These cubature rules appear to be new (see the list in \cite{CR}). Their numbers of nodes are  
more than the minimal given in \eqref{lwbd}. The existence of higher order Gauss-Lobatto rules
depends on the existence of minimal cubature rules of degree $2n-4$ for the weight $W_{1,1,1}$.
The latter cubature rules most likely do not exist for $n \ge 6$, but the algorithm can still
be applied to produce cubature formulas with $n-1$ points in each side of the triangle. 

It should be mentioned that there are cubature rules of degree $2n-1$ in the form of \eqref{cuba} 
that have fewer than $n-1$ points on each side, naturally with more interior points according to 
\eqref{lwbd2b}, such rules cannot be constructed directly by our algorithm. It is possible to modify
the algorithm, however, since $\CL_i$ for such rules cannot be positive definite on $\Pi_{n-1}$ but 
it must be positive definite on a subspace of $\Pi_{n-1}$. 

Finally, let us mention that our algorithm can also be modified for constructing cubature rules of degree 
$2n$. This means a cubature rule of degree $2n-3$ for $W_{\a+1,\b+1,\g+1}$ in Step 1, and 
a quadrature for $\CL_i$ of degree $2n-2$ in Step 2. The quadrature of degree $2n-2$ is generated 
by a quasi-orthogonal polynomial of the form $q_n := p_n + \a p_{n-1}$, where $\a$ is a free parameter
which can be fixed by requiring, say, $q_n(1/2) =0$, which means fixing the middle point on the 
corresponding side of the triangle as a node of the cubature. We have tried this construction for 
cubature rules of degree $6$ with 4 interior points, 3 on each sides and 1 at each vertex, starting with 
a cubature rule for $W_{1,1,1}$ of degree 3. The Lobatto type cubature rule of degree $6$ that we
obtained, however, has three negative weights. 

\medskip\noindent
{\bf Acknowledgment.} The author thanks both referees for their careful reading of the manuscript, 
especially David Day for his extensive and thoughtful suggestions.


\begin{thebibliography}{99}

\bibitem{BS}
        H. Berens and H. J. Schmid, 
        On the number of nodes of odd degree cubature formulae for integrals with Jacobi weight on a simplex,
        {\it Numerical Integration (Bergen, 1991)}, 37--44, NATO Adv. Sci. Inst. Ser. C Math. Phys. Sci., 357, 
        Kluwer Acad. Publ., Dordrecht, 1992.
      
\bibitem{BSX}
        H. Berens, H. Schmid and Y. Xu, 
        On two-dimensional definite orthogonal systems and on lower bound 
        for the number of associated cubature formulas,
        \textit{SIAM J. Math. Anal.} \textbf{26} (1995), 468--487.

\bibitem{CMS}
        R. Cools, I. P. Mysovskikh, and H. J. Schmid,
         Cubature formulae and orthogonal polynomials,
         {\it J. Comp. Appl, Math} {\bf 127} (2001) 121 - 152.

\bibitem{CR}
        R. Cools and P. Rabinowitz,
        Monomial cubature rules since ``Stroud'': a compilation. 
        \textit{J. Comput. Appl. Math.} \textbf{48} (1993), no. 3, 309--326. 

\bibitem{DX}
        C. F. Dunkl and Y. Xu,
         \textit{Orthogonal Polynomials of Several Variables}
         Encyclopedia of Mathematics and its Applications \textbf{81},
         Cambridge University Press, Cambridge, 2001.         

\bibitem{H} 
        B. T. Helenbrook, 
        On the existence of explicit $hp$-finite element method using Gauss-Lobatto
        integration on the triangle, 
         \textit{SIAM J. Numer. Anal.} \textbf{47} (2009), 1304-1318.
       
\bibitem{M}
        H. M\"oller,
        Kubaturformeln mit minimaler Knotenzahl,
	\textit{Numer. Math.} \textbf{ 25} (1976), 185--200.

\bibitem{St}
        A. H. Stroud, 
        {\it Approximate calculation of multiple integrals}, 
        Prentice-Hall, Inc., Englewood Cliffs, N.J., 1971.

\bibitem{X}
        Y. Xu,
        \textit{Common zeros of polynomials in several variables and higher     
        dimensional quadrature},
	Pitman Research Notes in Mathematics Series, Longman, Essex,
        1994.

\bibitem{X2}
        Y. Xu,
        On zeros of multivariate quasi-orthogonal polynomials and Gaussian cubature formulae, 
        {\it SIAM J. Math. Anal.} \textbf{25} (1994), 991-1001.

\bibitem{ZCL}
       L. Zhang, T. Cui, and H. Liu,
       A set of symmetric quadratures on triangles and tetrahedra,
       \textit{J. Computational Math.}, \textbf{27} (2009), 89 - 96. 

\end{thebibliography}
\end{document}